\magnification\magstep1
\newread\AUX\immediate\openin\AUX=\jobname.aux
\def\ref#1{\expandafter\edef\csname#1\endcsname}
\ifeof\AUX\immediate\write16{\jobname.aux gibt es nicht!}\else
\input \jobname.aux
\fi\immediate\closein\AUX
\def\today{\number\day.~\ifcase\month\or
  Januar\or Februar\or M{\"a}rz\or April\or Mai\or Juni\or
  Juli\or August\or September\or Oktober\or November\or Dezember\fi
  \space\number\year}
\font\sevenex=cmex7
\scriptfont3=\sevenex
\font\fiveex=cmex10 scaled 500
\scriptscriptfont3=\fiveex

\def\phi{\varphi}
\def\epsilon{\varepsilon}
\def\theta{\vartheta}
\def\uauf{\lower1.7pt\hbox to 3pt{%
\vbox{\offinterlineskip
\hbox{\vbox to 8.5pt{\leaders\vrule width0.2pt\vfill}%
\kern-.3pt\hbox{\lams\char"76}\kern-0.3pt%
$\raise1pt\hbox{\lams\char"76}$}}\hfil}}
\def\cite#1{\expandafter\ifx\csname#1\endcsname\relax
{\bf?}\immediate\write16{#1 ist nicht definiert!}\else\csname#1\endcsname\fi}
\def\expandwrite#1#2{\edef\next{\write#1{#2}}\next}
\def\neverexpand{\noexpand\noexpand\noexpand}
\def\strip#1\ {}
\def\ncite#1{\expandafter\ifx\csname#1\endcsname\relax
{\bf?}\immediate\write16{#1 ist nicht definiert!}\else
\expandafter\expandafter\expandafter\strip\csname#1\endcsname\fi}
\newwrite\AUX
\immediate\openout\AUX=\jobname.aux
\newcount\Abschnitt\Abschnitt0
\def\beginsection#1. #2 \par{\advance\Abschnitt1%
\vskip0pt plus.10\vsize\penalty-250
\vskip0pt plus-.10\vsize\bigskip\vskip\parskip
\edef\TEST{\number\Abschnitt}
\expandafter\ifx\csname#1\endcsname\TEST\relax\else
\immediate\write16{#1 hat sich geaendert!}\fi
\expandwrite\AUX{\neverexpand\ref{#1}{\TEST}}
\leftline{\bf\number\Abschnitt. \ignorespaces#2}%
\nobreak\smallskip\noindent\SATZ1}
\def\Proof:{\par\noindent{\it Proof:}}
\def\Remark:{\ifdim\lastskip<\medskipamount\removelastskip\medskip\fi
\noindent{\bf Remark:}}
\def\Remarks:{\ifdim\lastskip<\medskipamount\removelastskip\medskip\fi
\noindent{\bf Remarks:}}
\def\Definition:{\ifdim\lastskip<\medskipamount\removelastskip\medskip\fi
\noindent{\bf Definition:}}
\def\Example:{\ifdim\lastskip<\medskipamount\removelastskip\medskip\fi
\noindent{\bf Example:}}
\newcount\SATZ\SATZ1
\def\proclaim #1. #2\par{\ifdim\lastskip<\medskipamount\removelastskip
\medskip\fi
\noindent{\bf#1.\ }{\it#2}\Par
\ifdim\lastskip<\medskipamount\removelastskip\goodbreak\medskip\fi}
\def\Aussage#1{%
\expandafter\def\csname#1\endcsname##1.{\ifx?##1?\relax\else
\edef\TEST{#1\penalty10000\ \number\Abschnitt.\number\SATZ}
\expandafter\ifx\csname##1\endcsname\TEST\relax\else
\immediate\write16{##1 hat sich geaendert!}\fi
\expandwrite\AUX{\neverexpand\ref{##1}{\TEST}}\fi
\proclaim {\number\Abschnitt.\number\SATZ. #1\global\advance\SATZ1}.}}
\Aussage{Theorem}
\Aussage{Proposition}
\Aussage{Corollary}
\Aussage{Lemma}
\font\la=lasy10
\def\strich{\hbox{$\vcenter{\hbox
to 1pt{\leaders\hrule height -0,2pt depth 0,6pt\hfil}}$}}
\def\dashedrightarrow{\hbox{%
\hbox to 0,5cm{\leaders\hbox to 2pt{\hfil\strich\hfil}\hfil}%
\kern-2pt\hbox{\la\char\string"29}}}

\def\Bindestrich{\penalty10000-\hskip0pt}
\let\_=\Bindestrich
\def\.{{\sfcode`.=1000.}}

\def\Par{\par}
\def\:={\mathrel{\raise0,9pt\hbox{.}\kern-2,77779pt
\raise3pt\hbox{.}\kern-2,5pt=}}
\def\=:{\mathrel{=\kern-2,5pt\raise0,9pt\hbox{.}\kern-2,77779pt
\raise3pt\hbox{.}}}

\def\Ugleich{\hbox{$\cup$\kern.5pt\vrule depth -0.5pt}}
\def\|#1|{\mathop{\rm#1}\nolimits}
\def\<{\langle}
\def\>{\rangle}
\let\Times=\times
\def\times{\mathop{\Times}}
\let\Otimes=\otimes
\def\otimes{\mathop{\Otimes}}
\catcode`\@=11
\def\hex#1{\ifcase#1 0\or1\or2\or3\or4\or5\or6\or7\or8\or9\or A\or B\or
C\or D\or E\or F\else\message{Warnung: Setze hex#1=0}0\fi}
\def\fontdef#1:#2,#3,#4.{%
\alloc@8\fam\chardef\sixt@@n\FAM
\ifx!#2!\else\expandafter\font\csname text#1\endcsname=#2
\textfont\the\FAM=\csname text#1\endcsname\fi
\ifx!#3!\else\expandafter\font\csname script#1\endcsname=#3
\scriptfont\the\FAM=\csname script#1\endcsname\fi
\ifx!#4!\else\expandafter\font\csname scriptscript#1\endcsname=#4
\scriptscriptfont\the\FAM=\csname scriptscript#1\endcsname\fi
\expandafter\edef\csname #1\endcsname{\fam\the\FAM\csname text#1\endcsname}
\expandafter\edef\csname hex#1fam\endcsname{\hex\FAM}}
\catcode`\@=12 

\fontdef Ss:cmss10,,.
\fontdef Fr:eufm10,eufm7,eufm5.


%
\fontdef bbb:msbm10,msbm7,msbm5.
\fontdef mbf:cmmib10,cmmib7,.
\fontdef msa:msam10,msam7,msam5.

\def\NN{{\bbb N}}
\def\QQ{{\bbb Q}}

\def\ZZ{{\bbb Z}}

\def\cP{{\cal P}}
\def\cT{{\cal T}}

\mathchardef\leer=\string"0\hexbbbfam3F
\mathchardef\subsetneq=\string"3\hexbbbfam24
\mathchardef\semidir=\string"2\hexbbbfam6E
\mathchardef\dirsemi=\string"2\hexbbbfam6F
\mathchardef\haken=\string"2\hexmsafam78
\mathchardef\auf=\string"3\hexmsafam10
\let\OL=\overline
\def\overline#1{{\hskip1pt\OL{\hskip-1pt#1\hskip-1pt}\hskip1pt}}


%
\abovedisplayskip 9.0pt plus 3.0pt minus 3.0pt
\belowdisplayskip 9.0pt plus 3.0pt minus 3.0pt
\newdimen\Grenze\Grenze2\parindent\advance\Grenze1em
\newdimen\Breite
\newbox\DpBox
\def\NewDisplay#1$${\Breite\hsize\advance\Breite-\hangindent
\setbox\DpBox=\hbox{\hskip2\parindent$\displaystyle{#1}$}%
\ifnum\predisplaysize<\Grenze\abovedisplayskip\abovedisplayshortskip
\belowdisplayskip\belowdisplayshortskip\fi
\global\futurelet\nexttok\WEITER}
\def\WEITER{\ifx\nexttok\qed\expandafter\leftQEDdisplay
\else\leftdisplay\fi}
\def\leftdisplay{\hskip-\hangindent\leftline{\box\DpBox}$$}
\def\leftQEDdisplay{\hskip-\hangindent
\line{\copy\DpBox\hfill\lower\dp\DpBox\copy\QEDbox}%
\belowdisplayskip0pt$$\bigskip\let\nexttok=}
\everydisplay{\NewDisplay}
\newbox\QEDbox
\newbox\nichts\setbox\nichts=\vbox{}\wd\nichts=2mm\ht\nichts=2mm
\setbox\QEDbox=\hbox{\vrule\vbox{\hrule\copy\nichts\hrule}\vrule}
\def\qed{\leavevmode\unskip\hfil\null\nobreak\hfill\copy\QEDbox\medbreak}
\newdimen\HIindent
\newbox\HIbox
\def\setHI#1{\setbox\HIbox=\hbox{#1}\HIindent=\wd\HIbox}
\def\HI#1{\par\hangindent\HIindent\hangafter=0\noindent\leavevmode
\llap{\hbox to\HIindent{#1\hfil}}\ignorespaces}

\newdimen\maxSpalbr
\newdimen\altSpalbr

\def\beginrefs{%
\expandafter\ifx\csname Spaltenbreite\endcsname\relax
\def\Spaltenbreite{1cm}\immediate\write16{Spaltenbreite undefiniert!}\fi
\expandafter\altSpalbr\Spaltenbreite
\maxSpalbr0pt
\def\L|Abk:##1|Sig:##2|Au:##3|Tit:##4|Zs:##5|Bd:##6|S:##7|J:##8||{%
\edef\TEST{[##2]}
\expandafter\ifx\csname##1\endcsname\TEST\relax\else
\immediate\write16{##1 hat sich geaendert!}\fi
\expandwrite\AUX{\neverexpand\ref{##1}{\TEST}}
\setHI{[##2]\ }
\ifnum\HIindent>\maxSpalbr\maxSpalbr\HIindent\fi
\ifnum\HIindent<\altSpalbr\HIindent\altSpalbr\fi
\HI{[##2]}
\ifx-##3\relax\else{##3}: \fi
\ifx-##4\relax\else{##4}{\sfcode`.=3000.} \fi
\ifx-##5\relax\else{\it ##5\/} \fi
\ifx-##6\relax\else{\bf ##6} \fi
\ifx-##8\relax\else({##8})\fi
\ifx-##7\relax\else, {##7}\fi\Par}
\def\B|Abk:##1|Sig:##2|Au:##3|Tit:##4|Reihe:##5|Verlag:##6|Ort:##7|J:##8||{%
\edef\TEST{[##2]}
\expandafter\ifx\csname##1\endcsname\TEST\relax\else
\immediate\write16{##1 hat sich geaendert!}\fi
\expandwrite\AUX{\neverexpand\ref{##1}{\TEST}}
\setHI{[##2]\ }
\ifnum\HIindent>\maxSpalbr\maxSpalbr\HIindent\fi
\ifnum\HIindent<\altSpalbr\HIindent\altSpalbr\fi
\HI{[##2]}
\ifx-##3\relax\else{##3}: \fi
\ifx-##4\relax\else{##4}{\sfcode`.=3000.} \fi
\ifx-##5\relax\else{(##5)} \fi
\ifx-##7\relax\else{##7:} \fi
\ifx-##6\relax\else{##6}\fi
\ifx-##8\relax\else{ ##8}\fi\Par}
\def\Pr|Abk:##1|Sig:##2|Au:##3|Artikel:##4|Titel:##5|Hgr:##6|Reihe:{%
\edef\TEST{[##2]}
\expandafter\ifx\csname##1\endcsname\TEST\relax\else
\immediate\write16{##1 hat sich geaendert!}\fi
\expandwrite\AUX{\neverexpand\ref{##1}{\TEST}}
\setHI{[##2]\ }
\ifnum\HIindent>\maxSpalbr\maxSpalbr\HIindent\fi
\ifnum\HIindent<\altSpalbr\HIindent\altSpalbr\fi
\HI{[##2]}
\ifx-##3\relax\else{##3}: \fi
\ifx-##4\relax\else{##4}{\sfcode`.=3000.} \fi
\ifx-##5\relax\else{In: \it ##5}. \fi
\ifx-##6\relax\else{(##6)} \fi\PrII}
\def\PrII##1|Bd:##2|Verlag:##3|Ort:##4|S:##5|J:##6||{%
\ifx-##1\relax\else{##1} \fi
\ifx-##2\relax\else{\bf ##2}, \fi
\ifx-##4\relax\else{##4:} \fi
\ifx-##3\relax\else{##3} \fi
\ifx-##6\relax\else{##6}\fi
\ifx-##5\relax\else{, ##5}\fi\Par}
\bgroup
\baselineskip12pt
\parskip2.5pt plus 1pt
\hyphenation{Hei-del-berg}
\sfcode`.=1000
\beginsection References. References

}
\def\endrefs{%
\expandwrite\AUX{\neverexpand\ref{Spaltenbreite}{\the\maxSpalbr}}
\ifnum\maxSpalbr=\altSpalbr\relax\else
\immediate\write16{Spaltenbreite hat sich geaendert!}\fi
\egroup}


\def\rho{\varrho}

\fontdef Ss:cmss10,,.
\font\BF=cmbx10 scaled \magstep2
\font\CSC=cmcsc10 
\baselineskip15pt

{\baselineskip1.5\baselineskip\rightskip0pt plus 5truecm
\leavevmode\vskip0truecm\noindent
\BF A Recursion and a Combinatorial Formula for Jack Polynomials

}
\vskip1truecm
\leftline{{\CSC Friedrich Knop \& Siddhartha Sahi}%
\footnote*{\rm \rm The authors were partially supported by NSF grants.}}
\leftline{Department of Mathematics, Rutgers University, New Brunswick NJ
08903, USA}
\vskip1truecm
\beginsection Introduction. Introduction

The Jack polynomials $J_\lambda(x;\alpha)$ form a remarkable class of
symmetric polynomials. They are indexed by a partition $\lambda$ and
depend on a parameter $\alpha$. One of their properties is that several
classical families of symmetric functions can be obtained by
specializing $\alpha$, e.g., the monomial symmetric functions
$m_\lambda$ ($\alpha=\infty$), the elementary functions $e_{\lambda'}$
($\alpha=0)$, the Schur functions $s_\lambda$ ($\alpha=1$) and
finally the two classes of zonal polynomials ($\alpha=2$,
$\alpha=1/2$).

The Jack polynomials can be defined in various ways, e.g.:
\item{a)}as an orthogonal family of functions which is compatible with
the canonical filtration of the ring symmetric functions or
\item{b)}as simultaneous eigenfunctions of certain differential
operators (the Sekiguchi\_Debiard operators).

Recently Opdam, \cite{Op}, constructed a similar family
$F_\lambda(x;\alpha)$ of {\it non\_symmetric\/} polynomials. The index
runs now through all compositions $\lambda\in\NN^n$. They are defined
in a completely similar fashion, e.g., the Sekiguchi\_Debiard
operators are being replaced by the Cherednik differential\_reflection
operators (see section \cite{Operators}). It is becoming more and more
clear that these polynomials are as important as their symmetric
counterparts.

The purpose of this paper is to add to the existing characterizations of
Jack polynomials two further ones:
\item{c)}a recursion formula among the $F_\lambda$ together with
two formulas to obtain $J_\lambda$ from them.
\item{d)}combinatorial formulas of both $J_\lambda$ and $F_\lambda$ in
terms of certain generalized tableaux.

\noindent There are many advantages of these new characterizations over
the ones mentioned above. In a) and b), the
existence of functions with these properties is not obvious and
requires a proof whereas c) and d) could immediately serve as
a {\it definition\/} of Jack polynomials. Moreover, a) and b)
determine the functions only up to a scalar while c) and d) give
automatically the right normalization.

More importantly, our formulas are explicit enough such that both the recursion
relation and the combinatorial formula enable us to prove a conjecture of
Macdonald and Stanley (\cite{Mac}, \cite{St}). For a partition $\lambda$
let $m_i(\lambda)$ be the number of parts which are equal to $i$ and let
$u_\lambda:=\prod_{i\ge1}m_i(\lambda)!$. Then we prove

\Theorem. Let $J_\lambda(x;\alpha)=\sum_\mu
v_{\lambda\mu}(\alpha)m_\mu(x)$. Then all functions
$\tilde v_{\lambda\mu}(\alpha):=u_\mu^{-1}v_{\lambda\mu}(\alpha)$ are
polynomials in $\alpha$ with positive integral coefficients.

\noindent For an analogous statement for the $F_\lambda$ see
\cite{MacStConj}. We would like to mention the recent papers \cite{LV1}
and \cite{LV2} of Lapointe and Vinet which, by completely different
methods, establish that $v_{\lambda\mu}$ is a polynomial with integral
coefficients. Except for special cases, before that it is was not
even known that $v_{\lambda\mu}$ is a polynomial.

We continue with the description of c) and d).
First, the recursion formula.

For $\lambda\in\NN^n$ we define the {\it degree}
$|\lambda|:=\sum_i\lambda_i$. Its {\it length\/} $l(\lambda)$ is the
the maximal index $i$ such that $\lambda_i\ne0$. With
$m:=l(\lambda)$ we define $\tilde\lambda_m\:=\alpha\lambda_m+k+1$
where $k$ is the number of indices $i=1,\ldots,m-1$ with
$\lambda_i<\lambda_m$. Moreover, let
$\lambda^*:=(\lambda_m-1,\lambda_1,\ldots,\lambda_{m-1},0,\ldots,0)$.
For $i=m,\ldots,n$ let $$
f_i(x):=F_{\lambda^*}(x_i,x_1,\ldots,x_{i-1},x_{i+1},\ldots,x_n).
$$
Then we prove (\cite{Induc}):
$$
F_\lambda(x)=\tilde\lambda_m x_mf_m(x)+
x_{m+1}f_{m+1}(x)+x_{m+2}f_{m+2}(x)+\ldots+x_nf_n(x).
$$
The symmetric functions are most easily obtained if the number of
variables is big enough, i.e., $n\ge2m$. Let $\lambda^+\in\NN^{n-m}$
be the partition which is a permutation of
$(\lambda_1,\ldots,\lambda_{n-m})$. Then we prove (\cite{Sym2})
$$
J_{\lambda^+}(z_{m+1},\ldots,z_n)=
F_\lambda(0,\ldots,0,z_{m+1},\ldots,z_n).
$$

Now, we describe the combinatorial formula. For simplicity we
restrict ourselves to the symmetric case $J_\lambda$. Let
$\lambda$ be a partition. A {\it generalized tableau of shape
$\lambda$\it} is a labeling $T$ of the boxes in the Ferrers diagram of
$\lambda$ by numbers $1,2,\ldots,n$. To $T$, we associate the monomial
$x^T:=\prod_{s\in\lambda}x_{T(s)}$.

We call $T$ {\it admissible\/} if it satisfies for all
boxes $(i,j)\in\lambda$:
\item{a)} $T(i,j)\not= T(i',j)$ whenever $i'>i$
\item{b)} $T(i,j)\not= T(i',j-1)$ whenever $j>1$ and $i'<i$.

\noindent A box $s=(i,j)\in\lambda$ is {\it critical\/} (for $T$) if
$j>1$ and $T(i,j)=T(i,j-1)$.

\noindent Let $\lambda'$ be the dual partition to $\lambda$. The
armlength of $s=(i,j)\in\lambda$ is defined as
$a_\lambda(s):=\lambda_i-j$. Likewise, the leglength is defined as
$l_\lambda(s):=\lambda'_j-i$. Then we introduce the linear polynomial
$d_\lambda(s):=\alpha(a_\lambda(s)+1)+(l_\lambda(s)+1)$. With
$d_T(\alpha):=\prod_{s \|critical|}d_\lambda(s)$ our formula
reads (\cite{CombFor})
$$
J_\lambda(x;\alpha)=\sum_{T \|admissible|}d_T(\alpha)x^T.
$$

This formula immediately implies the Macdonald\_Stanley conjecture.
Consider a partition $\mu$ and the set $\cT$ of all tableaux $T$ with
$x^T=x^\mu$. Let $H$ be the group of permutations $\pi$ of the labels
$1,\ldots,n$ such that $\mu_{\pi(i)}=\mu_i$ for all $i$ and $\pi(i)=i$
whenever $\mu_i=0$. This group acts freely on $\cT$ by permuting the
labels such that $d_T(\alpha)$ and $x^T$ are left invariant. Since the
order of $H$ is $u_\mu$, we obtain that the coefficient of $x^\mu$ is
divisible by $u_\mu$.

In the sequel we prove first that the the eigenfunctions of the
Cherednik operators satisfy our recursion formula. Then we prove that
the functions defined by the combinatorial formula satisfy the recursion
relation as well.

\beginsection Jack polynomials. The definition of Jack polynomials

Most constructions and results in the following two sections can be
found in Opdam's paper~\cite{Op} in the framework of arbitrary root
systems. Here we are only interested in the case of ${\Ss
A}_{n-1}$.

Let $\cP\:=\QQ[x_1,\ldots,x_n]$ be the ring of polynomials. For an
indeterminate $\alpha$ let $\cP_\alpha=\cP\otimes_\QQ \QQ(\alpha)$.
If $\alpha$ is such that $1/\alpha$ is a non\_negative {\it integer}
then $\delta^{1/\alpha}(x)\:=\prod_{i\ne j}(1-x_ix_j^{-1})^{1/\alpha}$
is in the Laurent polynomial ring $\cP'=\cP[x^{-1}]$. Let $[f]_0\in\QQ$
denote the constant term of $f\in\cP'$. Then
$$
\<f,g\>_\alpha\:=[f(x)g(x^{-1})\delta^{1/\alpha}(x)]_0
$$
defines a non\_degenerate scalar product on $\cP$.

Consider $\Lambda\:=\NN^n$. The {\it degree} of
$\lambda=(\lambda_i)\in\Lambda$ is defined as
$|\lambda|\:=\sum_i\lambda_i$ and its {\it length} as
$l(\lambda)\:=\|max|\{k\mid\lambda(k)\ne0\}$ (with $l(0)\:=0$). We
recall the (partial) order relation of \cite{Op} on $\Lambda$. We start
with the usual ordering on the set $\Lambda^+\subseteq\Lambda$ of all
partitions $\lambda_1\ge\lambda_2\ge\ldots\ge\lambda_n\ge0$. Here
$\lambda\ge\mu$ if $|\lambda|=|\mu|$ and $$
\lambda_1+\lambda_2+\ldots+\lambda_i\ge\mu_1+\mu_2+\ldots+\mu_i
\quad\hbox{for all }i=1,\ldots,n.
$$
This order relation is extended to all of $\Lambda$ as follows.
Clearly, the symmetric group $W$ on $n$ letters acts on $\Lambda$ and
for every $\lambda\in\Lambda$ there is a unique partition $\lambda^+$
in the orbit $W\lambda$. For all permutations $w\in W$ with
$\lambda=w\lambda^+$ there is a unique one, denoted by $w_\lambda$, of
minimal length. We define $\lambda\ge\mu$ if either $\lambda^+>\mu^+$
or $\lambda^+=\mu^+$ and $w_\lambda\le w_\mu$ in the Bruhat order of
$W$. In particular, $\lambda^+$ is the unique {\it maximum\/} of
$W\lambda$.

Non\_symmetric Jack polynomials are defined by the following theorem.
Here $x^\lambda$ be the monomial  corresponding to $\lambda$.

\Theorem AAA. {\rm(\cite{Op}~2.6)} For every
$\lambda\in\Lambda$ there is a unique polynomial
$E_\lambda(x;\alpha)\in\cP_\alpha$ satisfying
\item{i)}$E_\lambda=x^\lambda+\sum_{\mu\in\Lambda:\mu<\lambda}
c_{\lambda\mu}(\alpha)x^\mu$ and
\item{ii)}$\<E_\lambda,x^\mu\>_\alpha=0$ for all $\mu\in\Lambda$ with
$\mu<\lambda$ and almost all $\alpha$ such that $1/\alpha\in\NN$.\Par
\noindent Moreover, the collection $\{E_\lambda\mid\lambda\in\Lambda\}$
forms a $\QQ(\alpha)$\_linear basis of $\cP_\alpha$.

The symmetric group $W$ acts on $\cP$ in the obvious way. Then $\cP^W$
is the algebra of symmetric functions. For $\lambda\in\Lambda^+$ let
$m_\lambda:=\sum Wx^\lambda$ denote the corresponding monomial
symmetric function. Then (symmetric) Jack polynomials are defined by:

\Theorem BBB. {\rm (\cite{Mac}~10.13)} For every $\lambda\in\Lambda^+$
there is a unique symmetric polynomial
$P_\lambda(x;\alpha)\in\cP_\alpha^W$ satisfying
\item{i)}$P_\lambda=m_\lambda+\sum_{\mu\in\Lambda^+:\mu<\lambda}
c'_{\lambda\mu}(\alpha)m_\mu$ and
\item{ii)}$\<P_\lambda,m_\mu\>_\alpha=0$ for all $\mu\in\Lambda^+$
with $\mu<\lambda$ and almost all $\alpha$ with $1/\alpha\in\NN$.\Par
\noindent
Moreover, the collection $\{P_\lambda\mid\lambda\in\Lambda\}$ forms a
$\QQ(\alpha)$\_linear basis of $\cP^W_\alpha$.

\noindent
An easy consequence of the definitions is:

\Lemma Space. For $\lambda\in\Lambda^+$ let
$\cP_\lambda\subset\cP_\alpha$ be the $\QQ(\alpha)$\_linear subspace
spanned by the $E_{w\lambda}$, $w\in W$. Then $\cP_\lambda$ is
$W$\_stable and $\cP_\lambda^W=\QQ(\alpha)P_\lambda$.

\noindent The action of $w\in W$ on $\cP_\lambda$ is, in general,
difficult to describe in terms of the basis $E_{w\lambda}$, but, for a
simple reflection $s_i\:=(i,i{+}1)\in W$, this is possible. 
We first present only a
special case and the rest later (\cite{LemA}).

\Lemma Symmetry. Let $\lambda\in\Lambda$ with
$\lambda_i=\lambda_{i+1}$. Then $s_iE_\lambda=E_\lambda$.

\Proof: This follows directly from the definition and the fact that
$\mu<\lambda=s_i\lambda$ implies $s_i\mu<\lambda$.\qed

\noindent One consequence of this lemma is that if $\lambda_i=0$ for
all $i>m$ then $E_\lambda$ is symmetric in the variables
$x_{m{+}1},\ldots,x_n$. This fact will be crucial later on.

\beginsection Operators. Definition of Cherednik's operators

As already mentioned, the symmetric group $W$ acts on $\cP$. For
$i\ne j$ let $s_{ij}\in W$ denote the transposition $(ij)$. Then
$$
N_{ij}:={1-s_{ij}\over x_i-x_j}
$$
is a well defined operator on $\cP$. Next, for $i=1,\ldots,n$ we
define the following differential\_reflection operators, which were
first studied by Cherednik~\cite{Che} (see also \cite{Op}):
$$
\xi_i:=\alpha x_i{\partial\over\partial x_i}+
\sum_{j=1}^{i-1}N_{ij}x_j+\sum_{j=i+1}^nx_jN_{ij}
$$
\Remark: The operators in \cite{Op} depend on the
choice of a positive root system. We use 
$\{-x_1+x_2,\ldots,-x_{n-1}+x_n\}$ as the set of simple roots. This has
the advantage that the $\xi_i$ are stable under adding variables.
\medskip\noindent
The $\xi_i$ commute pairwise. This is most easily seen by using
\cite{EigBas} below. Furthermore, they satisfy the following
commutation relations with the simple reflections
$s_i=s_{i\,i{+}1}$. This one checks by direct calculation.
$$
\eqalign{
\xi_is_i-s_i\xi_{i{+}1}&=1\cr
\xi_{i{+}1}s_i-s_i\xi_i&=-1\cr
\xi_is_j-s_j\xi_i&=0\qquad j\ne i,i+1\cr}
$$
(In other words, the $s_j$ and $\xi_i$ generate a graded Hecke
algebra.)

\Lemma EigVal. (a) The action of $\xi_i$ on $\cP$ is triangular. More precisely
$$
\xi_i(x^\lambda)=\bar\lambda_ix^\lambda+\sum_{\mu\in\Lambda:\mu<\lambda}
c_\mu x^\mu
$$
where $\bar\lambda_i\:=\alpha\lambda_i-(k_i'+k_i'')$ with
$$
\eqalign{
k_i'&=\#\{j=1,\ldots,i-1\mid\lambda_j\ge\lambda_i\}\cr
k_i''&=\#\{j=i+1,\ldots,n\mid\lambda_j>\lambda_i\}\cr}
$$
(b) For $1/\alpha\in\NN$, the operator $\xi_i$ is
symmetric with respect to the scalar product $\<\ ,\ \>_\alpha$.

\Proof: (a) is \cite{Op}~2.10 and (b) is \cite{Che}~3.8. The key to part
(a) is the observation that $(N_{ij}x_j)(x_i^ax_j^b)$ contains
$x_i^ax_j^b$ if and only if $a\le b$ while for
$(x_jN_{ij})(x_i^ax_j^b)$ one needs $a<b$.\qed

\Corollary EigBas. {\rm(\cite{Op}~2.7)} The $E_\lambda$ form a
simultaneous eigenbasis for the $\xi_i$. More precisely,
$\xi_i(E_\lambda)=\bar\lambda_iE_\lambda$.

\Remarks: 1. For an alternate proof for the existence of a simultaneous
eigenbasis see the remark after \cite{Induc} below.

\noindent 2. The eigenvalues $\bar\lambda_i$ could be more concisely
described as follows. Consider the vector $\rho\:=(0,-1,-2,\ldots,-n+1)$.
Then $\bar\lambda_i=(\alpha\lambda+w_\lambda\rho)_i$. 
\medskip
Another consequence is stability:

\Corollary Stable1. Let $\lambda\in\Lambda$ with $\lambda_n=0$ and
$\lambda'\:=(\lambda_1,\ldots,\lambda_{n{-}1})$. Then we have
$$
E_\lambda|_{x_n=0}=E_{\lambda'}\in \QQ(\alpha)[x_1,\ldots,x_{n{-}1}].
$$
If $\lambda$ is a partition, then
$$
P_\lambda|_{x_n=0}=P_{\lambda'}\in \QQ(\alpha)[x_1,\ldots,x_{n{-}1}].
$$

\Proof: Obviously, when substituting $x_n=0$, the operators
$\xi_1,\ldots,\xi_{n{-}1}$ induce their counterpart on
$\QQ(\alpha)[x_1,\ldots,x_{n{-}1}]$. Hence, the first statement follows
from \cite{EigBas} and then the second from \cite{Space}.\qed

\Remark: This Corollary allows to define $E_\lambda$ and
$P_\lambda$ in infinitely many variables $x_1$, $x_2$, $x_3,\ldots$
where $\lambda\in\NN^\infty$ is a sequence such that almost all
$\lambda_i$ are zero. More precisely, they lie in
$\cP^\infty\:=\|lim|\limits_{\longleftarrow}\QQ(\alpha)[x_1,\ldots,x_n]$
where the limit is to be taken in the category of graded algebras.
Actually, \cite{Symmetry} implies that the $E_\lambda$ even lie in the
subalgebra $\cP^{(\infty)}$ of {\it almost symmetric\/} functions,
i.e., those $f\in\cP^\infty$ which are symmetric in the variables
$x_m,x_{m{+}1},\ldots$ for some $m\ge1$ depending on $f$.

\beginsection Recursion. The recursion formula

We define ``creation operators'' for the $E_\lambda$. The first one is
very easy to define but seems to be new:
$$
\Phi:=x_ns_{n-1}s_{n-2}\ldots s_1,
$$
i.e.,
$$
(\Phi f)(x_1,\ldots,x_n)\:=x_n
f(x_n,x_1,\ldots,x_{n-1})\qquad(f\in\cP).
$$

\Lemma. The following relations hold:
$$
\eqalign{
\xi_i\Phi&=\Phi\xi_{i+1}\qquad\hbox{for }i=1,\ldots,n-1\cr
\xi_n\Phi&=\Phi(\xi_1+1)\cr}
$$

\Proof: Let $\tau=s_{n-1}\ldots s_1$. This is a cyclic
permutation with $x_n\tau=\tau x_1$. Then the assertion follows from
the following commutation relations which hold for all $1\le i\ne j<n$:
$$
\eqalign{
&x_i\partial_{x_i}\Phi=\Phi x_{i+1}\partial_{x_{i+1}},\quad 
x_n\partial_{x_n}\Phi=\Phi x_1\partial_{x_1}+\Phi.\cr
&N_{ij}x_j\Phi=\Phi N_{i{+}1\,j{+}1}x_{j{+}1},\quad
x_jN_{ij}\Phi=\Phi x_{j{+}1}N_{i{+}1\,j{+}1}\cr
&x_nN_{in}\Phi=x_nN_{in}x_n\tau=x_n\tau N_{i{+}1\,1}x_1=\Phi
N_{i{+}1\,1}x_1\cr
&N_{nj}x_j\Phi=N_{nj}x_nx_j\tau=
x_nx_jN_{nj}\tau=\Phi x_{j{+}1}N_{1\,{j{+}1}}\cr}
$$\qed

\Corollary Cyclic. Let $\lambda\in\Lambda$ with $\lambda_n\ne0$. Put
$\lambda^*\:=(\lambda_n-1,\lambda_1,\ldots,\lambda_{n-1})$. Then
$E_\lambda=\Phi(E_{\lambda^*})$.

\noindent Opdam \cite{Op}~1.2 constructed an operator which
permutes two entries:

\Proposition LemA. Let $i\in\{1,...,n-1\}$ and
$\lambda\in\Lambda$ with $\lambda_i>\lambda_{i+1}$. Then
$xE_\lambda=(xs_i+1)E_{s_i(\lambda)}$
with $x=\bar\lambda_i-\bar\lambda_{i+1}$.

\Proof: Let $E:=(xs_i+1)E_{s_i(\lambda)}$. Then one easily
verifies $\xi_j(E)=\bar\lambda_jE$ for all $j$. The assertion follows
by comparing the highest coefficient.\qed

\noindent
These operators together with $\Phi$ already suffice to generate all
$E_\lambda$, but we still have to divide by the factor $x$. We prove a
refinement.

\Lemma Creation. For $\lambda\in\Lambda$ with $1\le m\:=l(\lambda)\le n$ let
$\lambda^\sharp\:=(\lambda_1,\ldots,\lambda_{m-1},0,\ldots,0,\lambda_m)$.
Then
$(\bar\lambda_m+m)E_\lambda=X_\lambda(E_{\lambda^\sharp})$
where
$$
X_\lambda\:=(\bar\lambda_m+m)s_m\ldots s_{n-1}+
\sum_{i=m+1}^ns_is_{i+1}\ldots s_{n-1}
$$

\Proof: We prove the statement by induction on $n-m$, the number of
trailing zeros. If $m=n$ then $X_\lambda$ is just multiplication by
$(\bar\lambda_n+n)$. For $m=n-1$, the assertion follows from
\cite{LemA}. Assume now $m\le n-2$ and put $\lambda^\circ\:=
(\lambda_1,\ldots,\lambda_{m-1},0,\lambda_m,0\ldots,0)$. It follows
from \cite{EigVal} that $\bar\lambda_{m+1}=-m$ and
$\bar\lambda^\sharp_{m+1}=\bar\lambda_m$. Put $x\:=\bar\lambda_m+m$,
$\zeta_i=s_i\ldots s_{n-1}$, and
$\tau_j\:=\sum_{i=j+1}^n\zeta_i$. Then, by induction and \cite{LemA}, we
get $x(x+1)E_\lambda=(xs_m+1)
[(x+1)\zeta_{m+1}+\tau_{m+1}]E_{\lambda^\sharp}=
[x(x+1)\zeta_m+xs_m\tau_{m+1}+(x+1)\zeta_{m+1}+
\tau_{m+1}]E_{\lambda^\sharp}$.  Now we use that $s_m$
commutes with $\tau_{m+1}$ and that
$s_mE_{\lambda^\sharp}=E_{\lambda^\sharp}$ (\cite{Symmetry}). Thus
we obtain
$x(x+1)E_\lambda=(x+1)[x\zeta_m+\tau_{m+1}+\zeta_{m+1}]E_{\lambda^\sharp}=
(x+1)X_\lambda E_{\lambda^\sharp}$. Finally, $x+1\ne0$ since
$\lambda_m\ne0$.\qed

Now, we introduce another normalization of the Jack polynomials.
Recall that the {\it diagram\/} of $\lambda\in\Lambda$ is the set of
points (or {\it boxes})  $(i,j)\in\ZZ^2$ such that $1\le i\le n$ and
$1\le j\le\lambda_i$. As usual, we identify $\lambda$ with its diagram.
For each box $s=(i,j)\in\lambda$ we define the {\it arm\_length\/}
$a_\lambda(s)$, the {\it leg\_length\/} $l_\lambda(s)$ and the {\it
$\alpha$\_hooklengths\/} $c_\lambda(s)$, $d_\lambda(s)$ as follows:
$$
\eqalign{
a_\lambda(s)&\:=\lambda_i-j\cr
l'_\lambda(s)&\:=\#\{k=1,\ldots, i-1\mid j\le
\lambda_k+1\le\lambda_i\}\cr
l''_\lambda(s)&\:=\#\{k=i+1,\ldots,n\mid j\le
\lambda_k\le\lambda_i\}\cr
l_\lambda(s)&\:=l_\lambda'(s)+l_\lambda''(s)\cr
c_\lambda(s)&\:=\alpha a_\lambda(s)+(l_\lambda(s)+1)\cr
d_\lambda(s)&\:=\alpha(a_\lambda(s)+1)+(l_\lambda(s)+1)\cr}
$$
Now, we define
$$
\eqalign{
F_\lambda(x;\alpha)&
\:=\prod_{s\in\lambda}d_\lambda(s)E_\lambda(x;\alpha);\cr
J_\lambda(x;\alpha)&
\:=\prod_{s\in\lambda}c_\lambda(s)P_\lambda(x;\alpha).\cr}
$$
If $\lambda\in\Lambda^+$ is a partition then $l'(s)=0$ and
$l''(s)=l(s)$ is just the usual leg\_length. Moreover, $c_\lambda(s)$
is called the lower hook length in \cite{St}. This also shows that our
$J_\lambda(x;\alpha)$ coincides with $J_\lambda^{(\alpha)}$ in
\cite{Mac}.

First we state a simple lemma which calculates $d_\lambda(s)$ in a
special case.

\Lemma dPol. Let $\lambda\in\Lambda$ and $s=(i,1)\in\lambda$. Then
$d_\lambda(s)=\bar\lambda_i+i+a_0$ where
$a_0\:=\#\{k=i+1,\ldots,n\mid \lambda_k>0\}$.

\Proof: Follows directly from the definitions.\qed

Now we can prove our main recursion formula:

\Theorem Induc. For any $1\le k\le n$ put 
$$
\Phi_k\:=x_ks_{k-1}\ldots s_1
$$
For $\lambda\in\Lambda$ with $m:=l(\lambda)>0$ let
$$
\lambda^*\:=(\lambda_m-1,\lambda_1,\ldots,\lambda_{m-1},0,\ldots,0);
$$
$$
Y_\lambda\:=X_\lambda\Phi=
(\bar\lambda_m+m)\Phi_m+\Phi_{m+1}+\ldots+\Phi_n
$$
Then $F_\lambda=Y_\lambda(F_{\lambda^*})$.

\Proof: \cite{Cyclic} and \cite{Creation} imply
$xE_\lambda=Y_\lambda(E_{\lambda^*})$ with $x=\bar\lambda_m+m$.
The diagram of $\lambda^*$ is obtained from $\lambda$ by
taking the last non\_empty row, removing its first box $s_0=(m,1)$ and
putting the rest on top. One easily checks from the definitions that
the arm\_length and the leg\_length of the remaining boxes
don't change. Moreover $x=d_\lambda(s_0)$ by \cite{dPol}. This proves
the theorem.\qed

\Remark: One could use \cite{Induc} as a {\it definition\/} of
$F_\lambda$. Then reading the proofs of \cite{Creation} and
\cite{Induc} backwards one sees that the so defined functions are
simultaneous eigenfunctions for the Cherednik operators. This gives an
alternate proof of \cite{EigBas} and of the commutativity of Cherednik
operators.
\medskip

The following Corollary shows another way to normalize non\_symmetric
Jack polynomials in the case the number of variables is large enough.
It is an analogue of Stanley's normalization in \cite{St}.

\Corollary CoeffNS. Let $\lambda\in\Lambda$, put $d\:=|\lambda|$ and
$m\:=l(\lambda)$. Assume $n\ge m+d$. Then the coefficient of
$x_{m{+}1}\ldots x_{m{+}d}$ in $F_\lambda$ is $d!$.

\Proof: We have
$$
c:=x_{m{+}1}\ldots x_{m{+}d}=\Phi_i(x_{m+2},\ldots,x_{m+d})
$$
for $i=m+1,\ldots,m+d$ and this is the only way, $c$ can arise as the
image of an operator $\Phi_i$. Hence, \cite{Induc} implies that the
coefficient of $c$ in $F_\lambda$ is $d$ times the coefficient of
$x_{m+2}\ldots x_{m+d}$ in $F_{\lambda^*}$. But $F_{\lambda^*}$ is
symmetric in the variables $x_{m^*+1},\ldots,x_n$ where
$m^*=l(\lambda^*)$. The assertion follows by induction.\qed

We give the first of two ways how to obtain the symmetric Jack
polynomials from the non\_symmetric ones. Before we do so, recall some
notation. For any $\lambda\in\Lambda$ let
$m_i(\lambda)\:=\#\{k\mid\lambda_k=i\}$ and
$u_\lambda\:=\prod_{i\ge1}m_i(\lambda)!$.

\Theorem Sym1. For $\lambda\in\Lambda^+$ with $m:=l(\lambda)$ put
$$
\lambda^0\:=(\lambda_m-1,\ldots,\lambda_1-1,0,\ldots,0).
$$
Then $J_\lambda(x;\alpha)={1\over(n-m)!}\sum_{w\in
W}w\Phi^m(F_{\lambda^0})$.

\Proof: Denote the right hand side by $J$. Let
$\lambda^-\:=(0,\ldots,0,\lambda_m,\ldots,\lambda_1)$. Then
\cite{Cyclic} implies that $F':=\Phi^m(F_{\lambda^0})$ is proportional to
$F_{\lambda^-}$. \cite{Space} implies that $J$ and $J_\lambda$ are
proportional.

To see that they are equal, it suffices to compare the coefficients of
$x^{\lambda^-}$. Since $\Phi$ does not change the leading coefficient,
the coefficient of $z^{\lambda^-}$ in $F'$ is
$\prod_{s\in\lambda^-\atop s\ne(i,1)}d_{\lambda^-}(s)$. 
Since $\lambda^-$ is minimal in
$W\lambda$, no other monomial occuring in $F'$ is conjugated to
$x^{\lambda^-}$. Moreover, $F'$ is invariant for the isotropy group
$W_{\lambda^-}$. Its order is $(n-m)!u_\lambda$. Hence the
coefficient of $x^{\lambda^-}$ in $J$ is
$$
u_\lambda\prod_{s\in\lambda^-\atop s\ne(i,1)}d_{\lambda^-}(s)
$$
On the other hand, by definition, the coefficient of $x^{\lambda^-}$ in
$J_\lambda$ is
$$
\prod_{s\in\lambda}c_\lambda(s).
$$
Let $w\in W$ be the shortest permutation with $w(\lambda)=\lambda^-$.
This means $w(i)>w(j)$ whenever $\lambda_i>\lambda_j$ but
$w(i)<w(j)$ for $\lambda_i=\lambda_j$ and $i<j$. Consider the 
following correspondence between boxes of $\lambda$ and $\lambda^-$:
$$
\lambda\owns s=(i,j)\leftrightarrow s^-=(\pi(i),j+1)\in\lambda^- \ .
$$
This is defined for all $s$ with $j<\lambda_i$. One easily verifies
that $a_\lambda(s)=a_{\lambda^-}(s^-)+1$ and
$l_\lambda(s)=l_{\lambda^-}(s^-)$. Hence,
$c_\lambda(s)=d_{\lambda^-}(s^-)$, i.e., $s$ and $s^-$ contribute the
same factor to the products above. What is left out of the correspondence
are those boxes of $\lambda$ with $j=\lambda_i$ and the first column of
$\lambda^-$. The first type of these boxes contributes $u_\lambda$ to
the factor of $J_\lambda$. The second type doesn't  contribute by
construction. This shows that $J_\lambda=J$.\qed

\noindent This proof gives a bit more, namely a result of Stanley
(\cite{St}~Thm.~1.1 in conjunction with Thm.~5.6).

\Corollary CoeffS. Let $\lambda\in\Lambda^+$ with
$d:=|\lambda|\le n$. Then the coefficient of $m_{1^d}$ in
$J_\lambda$ is $d!$.

\Proof: We keep the notation of the proof of \cite{Sym1}. We have
$$
F'=\Phi^m(F_{\lambda^0})=x_{n{-}m{+}1}\ldots x_n
F_{\lambda_0}(x_{n{-}m{+}1},\ldots,x_n,x_1,\ldots,x_{n-m}).
$$
Hence every monomial occuring in $F'$ which contains each variable with
a power of at most one is of the form $x_{i_1}\ldots
x_{i_{d-m}}x_{n{-}m{+}1}\ldots x_n$ with $1\le i_1<\ldots<i_{d-m}\le
n-m$. By \cite{CoeffNS}, each of them has the coefficient $(d-m)!$. Hence
the coefficient of $x_1\ldots x_d$ in $J_\lambda$ is
$$
{1\over(n-m)!}(d-m)!{n-m\choose d-m}d!(n-d)!=d!.
$$\qed

\noindent The next theorem establishes a direct relation 
between symmetric and non\_symmetric Jack polynomials. It needs more
variables than symmetrization but has the advantage of being stable in
$n$. Observe, that $\lambda$ is not required to be a partition.

\Theorem Sym2. Let $\lambda\in\Lambda$ and $m\in\NN$ with $l(\lambda)\le
m\le n-l(\lambda)$. Let $\lambda^+$ be the unique partition which is a
permutation of $(\lambda_1,\ldots,\lambda_{n-m})$. Then
$$
J_{\lambda^+}(x_{m+1},\ldots,x_n)=
F_\lambda(0,\ldots,0,x_{m+1},\ldots,x_n).
$$

\Proof: Recall that $\cP_\lambda\subset\cP_\alpha$ is the
$\QQ(\alpha)$\_linear subspace spanned by all $E_{w\lambda}$, $w\in W$.
Then \cite{Stable1} implies that
$\cP_\lambda|_{x_{n{+}1{-}m}=\ldots=x_n=0}=\cP_{\lambda^+}\subseteq
\QQ(\alpha)[x_1,\ldots,x_{n{-}m}]$. Since $\cP_\lambda$ is $W$\_stable
we conclude that also
$\cP_\lambda|_{x_1=\ldots=x_m=0}=\cP_{\lambda^+}\subseteq
\QQ(\alpha)[x_{m{+}1},\ldots,x_n]$. \cite{Space} implies that both sides
of the equation are equal up to a factor $c\in\QQ(\alpha)$. To
determine $c$ we may assume that $n\ge m+|\lambda|$. Then, by
Corollaries \ncite{CoeffNS} and \ncite{CoeffS}, the monomial
$x_{m{+}1}\ldots x_{m+|\lambda|}$ figures on both sides with the same
non\_zero coefficient. Hence $c=1$.\qed

Although, as already indicated in the introduction, the
Macdonald\_Stanley conjecture follows immediately from the combinatorial
formula of the next section, a direct proof using the recursion formula
might be of interest. To formulate its analogue for the non\_symmetric
polynomials we introduce the following notation. Fix an $m\in\NN$ with
$0\le m\le n$. We split every $\lambda\in\Lambda$ in two parts
$\lambda'$ and $\lambda''$ where $\lambda'$ (respectively $\lambda''$)
consists of the first $m$ (respectively last $n-m$) components of
$\lambda$. We write $\lambda=\lambda'\lambda''$. Then we define the
partially symmetric monomial functions as $m_\lambda^{(m)}\:=\sum_\mu
x^{\lambda'\mu}$ where $\mu$ runs through all permutations of
$\lambda''$. Their augmented version is $\tilde
m_\lambda^{(m)}\:=u_{\lambda''}m_\lambda^{(m)}$. Let
$\Lambda^{(m)}\subseteq\Lambda$ be the set of those $\lambda$ where
$\lambda''$ is a partition. Observe that $\Lambda^{(0)}=\Lambda^+$,
$m_\lambda^{(0)}=m_\lambda$, and $\tilde m_\lambda^{(0)}=\tilde
m_\lambda$.

\Theorem MacStConj. a) Let $\lambda\in\Lambda$ and $m\in\NN$ with
$m\ge l(\lambda)$. Then
$$
F_\lambda(x;\alpha)=
\sum_{\mu\in\Lambda^{(m)}}a_{\lambda\mu}(\alpha)\tilde m_\mu^{(m)}
$$
with $a_{\lambda\mu}\in\NN[\alpha]$ for all $\mu\in\Lambda^{(m)}$.\Par
b) Let $\lambda\in\Lambda^+$. Then $J_\lambda(x;\alpha)=
\sum_{\mu\in\Lambda^+}b_{\lambda\mu}(\alpha)\tilde m_\mu$ with
$b_{\lambda\mu}\in\NN[\alpha]$ for all $\mu\in\Lambda^+$.

\Proof: Part b) follows immediately from a) and \cite{Sym2}. The proof of
a) is by induction on
$|\lambda|$. First observe that it suffices to prove the theorem for
$m=l(\lambda)$. Since $|\lambda^*|=|\lambda|-1$ and $l(\lambda^*)\le
m$, the assertion is true for $F_{\lambda^*}$. With
$\Psi\:=\Phi_{m+1}+\ldots+\Phi_n$ we have
$Y_\lambda=(\bar\lambda_m+m)\Phi_m+\Psi$. Moreover
$\bar\lambda_m+m=\alpha\lambda_i-k+m$ where $k$ is the number of
$j=1,\ldots,m-1$ with $\lambda_j\ge\lambda_m$. Thus $-k+m\ge1$.
By the recursion formula \cite{Induc}, it suffices to prove the following

\medskip
\noindent{\it Claim:\/} Let $\mu\in\Lambda^{(m)}$. Then both
$\Phi_m(\tilde m_\mu^{(m)})$ and $\Psi(\tilde m_\mu^{(m)})$ are linear
combinations of $\tilde m_\nu^{(m)}$, $\nu\in\Lambda^{(m)}$
with coefficients in $\NN$.

\medskip\noindent The effect of $\Phi_i$ on
monomials is $\Phi_i(x^\nu)=x^{\bar\nu}$ where
$$
\bar\nu\:=(\nu_2,\ldots,\nu_i,\nu_1+1,\nu_{i+1},\ldots,\nu_n).
$$
In particular, $\Phi_m$ affects only the first $m$ variables which
proves the claim for $\Phi_m$.

It is easy to check that the $\Phi_i$ satisfy the following
commutation relations:
$$
\vbox{\halign{$#$\hfill$\,=$&$\,#$,\hfill\quad&if $#$\hfill\cr
s_j\Phi_i&\Phi_is_j&i<j\cr
s_j\Phi_j&\Phi_{j+1}&\omit\cr
s_j\Phi_{j+1}&\Phi_j&\omit\cr
s_j\Phi_i&\Phi_is_{j+1}&i>j+1\cr}}
$$
This shows that $\Psi(\tilde m_\mu^{(m)})$ is invariant for
$s_{m{+}1},\ldots,s_n$. In particular, it suffices to check
the coefficient of $x^{\bar\nu}$ in $\Psi(\tilde m_\mu^{(m)})$
when $\bar\nu\in\Lambda^{(m)}$.

Assume that $x^\nu$ occurs in
$m_\mu^{(m)}$, i.e., that $\nu'=\mu'$ and that $\nu''$ is a
permutation of $\mu''$. Then $\nu$ is recovered from $\bar\nu$ by
removing a part $\bar\nu_i$ of $\bar\nu$ with $\bar\nu_i=k:=\mu_1+1$
and $i\ge m$ and putting $\mu_1$ in front. This shows that
$\Psi(m_\mu^{(m)})$ contains $x^{\bar\nu}$ with multiplicity
$m_k(\bar\nu'')$. Furthermore, $m_i(\bar\nu'')\le m_i(\mu'')$ for $i\ne
k$ and $m_k(\bar\nu'')\le m_k(\mu'')+1$. This implies that
$\Psi(\tilde m_\mu^{(m)})$ contains $\tilde m_\nu^{(m)}$ with positive
integral multiplicity.\qed

\def \lam {\lambda}

\beginsection Combinatorial. The combinatorial formula

In this section we give a simple and explicit formula for both the
symmetric and non\_symmetric Jack polynomials. Let $\Lambda\in\Lambda$.
A {\it generalized tableau\/} of shape $\lambda$ is a labeling $T$ of
the diagram of $\lambda$ by the numbers $1,\ldots,n$. The {\it weight\/}
of $T$ is $|T| =(|T|_1,\cdots, |T|_n)$ where $|T|_i$  is the number of
occurrences of the label $i$ in $T$. Of course $|T|$ is $S_n$-conjugate
to a unique partition. One writes $x^T$ for the monomial $x^{|T|}$.

\Definition: A generalized tableau of shape $\lam\in\Lambda$ is 
{\it admissible\/} if for all $(i,j)\in\lam$
\item{a)} $T(i,j)\ne T(i',j)$ if $i'>i.$
\item{b)} $T(i,j)\ne T(i',j-1)$ if $j>1$, $i'<i$.

\noindent
It is called {\it $0$\_admissible\/} if additionally
\item{c)} $T(i,j)\in \{i,i+1,\cdots,n\}$ if $j=1$.

\Definition: Let $T$ be a generalized tableau of shape $\lam$.
\item{a)} A point $(i,j)\in\lam$ is called {\it critical\/} if $j>1$ and
$T(i,j)=T(i,j-1)$.
\item{b)} The point $(i,j)\in\lam$ is called {\it $0$\_critical\/} if it is
critical or $j=1$ and $T(i,j)=i$.
\item{}The {\it hook-polynomials\/} of $T$ are
$$
\eqalign{
d_T(\alpha)&:=\prod_{s \|critical|}d_\lam(s,\alpha);\cr
d_T^0(\alpha)&:=\prod_{s\ 0\|-critical|}d_\lam(s,\alpha).\cr}
$$

\noindent
Our terminology can be explained as follows. Consider the tableau $T^0$
which arises from $T$ by adding a zero-th column and labeling its boxes
consecutively by $1,2,\ldots,n$. Then $T$ is $0$\_admissible if $T^0$
is admissible and a box $s$ in $T$ is $0$\_critical if it is critical in
$T^0$.

Our main theorem is:

\Theorem CombFor. Let $\lambda\in\Lambda$. Then
$$
F_\lam(x;\alpha)=\sum_{T 0\|-admissible|}d_T^0(\alpha) x^T.
$$
Let $\lambda^+\in\Lambda^+$ be the unique partition conjugated to
$\lambda$. Then
$$
J_{\lam^+}(x;\alpha)=\sum_{T \|admissible|}d_T(\alpha) x^T.
$$

\Proof: We prove first the formula for $J_{\lambda^+}$ assuming it for
$F_\lambda$. Let $m:=l(\lambda)$ and assume
$n\ge|\lambda|+l(\lambda)$. Consider only those tableaux of shape
$\lambda$ which contain only labels $>m$. Then ``$0$\_admissible'',
``0\_critical'' are the same as ``admissible'', ``critical'' respectively.
By \cite{Sym2}, the formula for $F_\lambda$ implies that for
$J_{\lambda^+}$.

For the non\_symmetric case, denote the right hand side of the formula
by $F'_\lambda$. We are going to prove the following two lemmas. 

\Lemma L2. Suppose $\lam_i=0$ and $\lam_{i+1}>0$, and write 
$d:=d_\lam(\alpha,(i+1,1))$. Then we have 
$d F'_{s_i\lam}=(d-1)s_i(F'_\lam) +F'_\lam$.

\noindent For $\lam\in\Lambda$  let $\Phi(\lam):=
(\lam_2,\cdots,\lam_n,\lam_1+1)$.

\Lemma L3. Let $d:= d_{\Phi\lam}(\alpha,(n,1))$ then $F'_{\Phi\lam}=
d \Phi(F'_\lam)$.

\noindent We finish first the proof of \cite{CombFor}. In the situation
of \cite{L2} let $\mu:=s_i\lam$. Then $d_\mu(i,1)=d-1$ while the
hook\_length of the remaining boxes doesn't change. Hence, if
$F_\lambda=c E_\lambda$ then $F_\mu={d-1\over d}cE_\mu$. Let $a_0$ be
the number of $k=i+2,\ldots, n$ with $\lambda_k>0$. Then
$\bar\mu_{i+1}=-i-a_0$ while $\bar\mu_i=d-1-i-a_0$ (\cite{dPol}).
Hence, $x:=\bar\mu_i-\bar\mu_{i+1}=d-1$.  With \cite{LemA} we get
$dF_\mu=(d-1)cE_\mu=xcE_\mu=(xs_i+1)cE_\lam=((d-1)s_i+1)F_\lam$.
We conclude from \cite{L2} that $F_\lam=F'_\lam$
implies $F_{s_i\lam}=F'_{s_i\lam}$.

In the same manner, we obtain from \cite{L3} and \cite{Cyclic} that 
$F_{\Phi\lam}=F'_{\Phi\lam}$ if and only if $F_\lam=F'_\lam$.
Since every $\lambda\in\Lambda$ is obtained by repeatedly applying
$\Phi$ or switching a zero and a non\_zero entry, the theorem follows by
induction (and $F_0=F'_0=1$).\qed

\noindent {\it Proof of \cite{L2}}: If $T$ is a tableau of shape $\lam$,
let $T'$ be the  tableau of shape $s_i\lam$ obtained by moving all the
points in row  $i+1$ up one unit to the previously empty row $i$.

Let us ignore for a moment the labels of $(i+1,1)\in\lam$ and 
$(i,1)\in s_i\lam$. For all {\it other\/} points in $T$, the
label is admissible (resp. critical) if and only it is so
for the corresponding point in $T'$, and the twisted hooklengths
are unchanged. (In  fact, $l'$, $l''$ and $a_\lam$ are all unchanged!)

To examine the contributions of $(i,1)$ and $(i+1,1)$,
we divide admissible tableaux $T$ of shape $\lam$ into two classes:

\noindent
$A= \{T\mid T(i+1,1)\ne i+1\}$, and $B=\{T\mid T(i+1,1)= i+1\}.$

Similarly we divide admissible tableaux $U$ of shape $s_i\lam$ into 
three classes:

\noindent
$A'= \{U\mid U(i,1)\ne i,i+1\}$,
$B'= \{U\mid U(i,1)=i+1\}$, and 
$B''= \{U\mid U(i,1)=i\}.$

The map $T \mapsto T'$ is a bijection from $A$ to $A'$,
and satisfies $d_T(\alpha)x^T=d_{T'}(\alpha)x^{T'}$. Also, if $T\in A$
then  replacing each occurrence of the label $i$ by $i+1$ and vice versa,
we get another tableau $s_iT\in A$ with $d_T(\alpha)=d_{s_iT}(\alpha)$.
This implies
$$\sum_{U\in A'}d_U(\alpha)x^U=
\sum_{T\in A}d_T(\alpha)x^T= s_i\sum_{T\in A}d_T(\alpha)x^T.$$

$T\mapsto T'$ is also a bijection from $B$ to $B'$,
however $T(i+1,1)$ is critical but $T'(i,1)$ is not. 
Since $d_\lam(\alpha,(i+1,1))=d$, we get 
$$d \sum_{U\in B'}d_U(\alpha)x^U=
 \sum_{T\in B}d_T(\alpha)x^T.$$

Finally $T\mapsto s_iT'$ is a bijection from $B$ to $B''$,
and $T(i+1,1)$, $s_iT'(i,1)$ are both critical. Since
$d_{s_i\lam}(\alpha,(i,1))=d-1$
($l''$, $a_\lam$ are unchanged, while $l'$ decreases by $1$), we get
$$d \sum_{U\in B''}d_U(\alpha)x^U=
(d-1) s_i\sum_{T\in B}d_T(\alpha)x^T.$$

Combining these we get
$dF'_{s_i\lam}=d\sum_{A'}+d\sum_{B'}+d\sum_{B''}=
 [\sum_A +(d-1)s_i\sum_A ]+ \sum_B +(d-1)s_i\sum_B=
(d-1)s_iF'_\lam+ F'_\lam.$ \qed

\noindent{\it Proof of \cite{L3}}:  For a tableau $T$ of shape $\lam$,
let $T'$ be the tableau of shape $\Phi\lam$ constructed 
as follows: 
\item{ 1)} move rows $2$ through $n$ up one place.
\item{ 2)} prefix the first row by a point with the label $1$
and move the row to the $n$-th place.
\item{ 3)} modify the labels by changing all $1$'s to $n$'s
and the other $i$'s to $(i-1)$'s.

If $s$ is a point in $T$, we write $s'$ for the corresponding 
point in $T'$, thus $s=(1,j)$ corresponds to $s'=(n,j+1)$ and
for $i>1$, $s=(i,j)$ corresponds to $s'=(i-1,j)$.

First observe the twisted hooklengths of corresponding points 
are the same. Indeed $a_\lam(s)= a_{\Phi\lam}(s')$,
and $l_\lam'(s)+l_\lam''(s)= l_{\Phi\lam}'(s')+l_{\Phi\lam}''(s')$.
($l'$ might decrease by $1$, but then $l''$ increases by $1$,
so that the sum is unchanged.) 

Second, note that if $T$ is admissible then so is $T'$. This is 
obvious for the first column, and for $(i,j)$ with $j>1$ and 
$i<n$, we only need to check that $T'(i,j)\ne T'(n,j))$. But 
these labels are obtained by applying 3) to the labels 
$T(i+1,j)$ and $T(1,j-1)$ which are distinct by the admissibility 
of $T$. The argument for the admissibility of $T'(n,j)$ is similar.

Next, note that the map $T\mapsto T'$ is actually a bijection from 
admissible tableaux of shape $\lam$ to those of shape $\Phi\lam$. 
The inverse map is obtained by deleting the label $T'(n,1)$
(which must be $n$), moving the last row to the top, and applying the 
inverse of 3).

Now, observe that the point $(n,1)$ is a critical point 
of $T'$, and any other point $s'$ of $T'$ is critical if and 
only if the corresponding point $s$ in $T$ is critical. This
is obvious for all points except $(n,2)$ which corresponds 
to $(1,1)$ in $T$; but $T'(n,2)=T'(n,1)=n$ if and only 
if $T(1,1)=1$.

Finally by 2) and 3), if the weight of $T$ is $\mu$ 
then the weight of $T'$ is $\Phi\mu$, thus $x^{T'}=\Phi(x^T)$.
This means $F'_{\Phi\lam}=\sum_{T'}d_{T'}(\alpha)x^{T'} = 
d_{\Phi\lam}(\alpha,(n,1))\sum_{T}d_T(\alpha)\Phi(x^T)
=d\Phi F'_\lam.$ \qed

\beginrefs

\L|Abk:Che|Sig:C|Au:Cherednik, I.|Tit:A unification of the
Knizhnik\_Zamolodchikov equations and Dunkl operators via
affine Hecke algebras|Zs:Invent. Math.|Bd:106|S:411--432|J:1991||

\L|Abk:LV1|Sig:LV1|Au:Lapointe, Luc; Vinet, Luc|Tit:Exact operator
solution of the Calogero\_Sutherland model%
|Zs:CRM-Preprint|Bd:2272|S:35 pages|J:1995||

\L|Abk:LV2|Sig:LV2|Au:Lapointe, Luc; Vinet, Luc|Tit:A Rodrigues formula
for the Jack polynomials and the Macdonald\_Stanley conjecture%
|Zs:CRM-Preprint|Bd:2294|S:5 pages|J:1995||

\B|Abk:Mac|Sig:M|Au:Macdonald, I.|Tit:Symmetric functions and Hall
polynomials (2nd ed.)|Reihe:-|Verlag:Cla\-rendon Press|Ort:Oxford|J:1995||

\L|Abk:Op|Sig:O|Au:Opdam, E.|Tit:Harmonic analysis for certain
representations of graded Hecke algebras%
|Zs:Acta Math.|Bd:175|S:75--121|J:1995||

\L|Abk:St|Sig:S|Au:Stanley, R.|Tit:Some combinatorial properties of
Jack symmetric functions|Zs:Advances Math.|Bd:77|S:76--115|J:1989||

\endrefs
\bye